\documentclass[preprint]{imsart}
\RequirePackage[OT1]{fontenc}
\RequirePackage{amsthm,amsmath}
\RequirePackage[numbers]{natbib}
\RequirePackage[colorlinks,citecolor=blue,urlcolor=blue]{hyperref}

\arxiv{}

\startlocaldefs
\numberwithin{equation}{section}
\theoremstyle{plain}
\newtheorem{thm}{Theorem}[section]
\newtheorem{lem}{Lemma}[section]
\newtheorem{defn}{Definition}[section]
\newtheorem{example}{Example}[section]
\newtheorem{rem}{Remark}[section]
\newtheorem{proposition}{Proposition}[section]
\newtheorem{fact}{Fact}[section]

\newtheorem{cor}{Corollary}[section]

\endlocaldefs

\begin{document}

\begin{frontmatter}
\title{Representation of the Dirac delta function in $\mathcal{C}(R^{\infty})$ in terms of infinite-dimensional Lebesgue measures  }
\runtitle{Representation of the Dirac delta function in $\mathcal{C}(R^{\infty})$ }

\begin{aug}
\author{\fnms{Gogi} \snm{Pantsulaia}\thanksref{t1}\ead[label=e1]{g.pantsulaia@gtu.ge}}
\and
\author{\fnms{Givi } \snm{Giorgadze}\ead[label=e2]{g.giorgadze}}

\affiliation{}

\thankstext{t1}{The research for this paper was partially supported by Shota
Rustaveli National Science Foundation's Grant no FR/116/5-100/14}
\runauthor{G.Pantsulaia and G.Giorgadze}

\address{I.Vekua Institute of Applied Mathematics, Tbilisi - 0143, Georgian Republic \\
\printead{e1}\\
Georgian Technical University, Tbilisi - 0175, Georgian Republic\\
  \printead*{e2}
}
\end{aug}

\begin{abstract}A representation of  the Dirac delta function in $ \mathcal{C}(R^{\infty})$ in terms of infinite-dimensional Lebesgue measures in $R^{\infty}$ is obtained and some it's properties are studied in this paper.
\end{abstract}



\end{frontmatter}

\section{Introduction}

 The Dirac delta function($\delta$-function) was introduced by Paul Dirac at the end of the 1920s in an effort to create the mathematical tools for the development of quantum  field theory. He referred to it as an improper functional � in Dirac (1930). Later, in 1947, Laurent Schwartz gave it a more rigorous mathematical definition as a spatial linear functional on the space of test functions $D$ (the set of all real-valued infinitely differentiable functions with compact support). Since the delta function is not really a function in the classical sense, one should not consider the value of the delta function at $x$. Hence, the domain of the delta function is $D$ and its value for $f \in D$ is $f(0)$. Khuri (2004) studied some interesting applications of the delta function in statistics.

 The purpose of the present paper is an introduction of a  concept of  the Dirac delta function in the class of all continuous functions defined in the infinite-dimensional topological vector space of all real valued sequences $R^{\infty}$ equipped with Tychonoff topology and a representation of this functional  in terms of  infinite-dimensional Lebesgue measures in $R^{\infty}$.

The paper is organized as follows:

In Section 2 we present a concept of ordinary and standard Lebesgue measures in $R^{\infty}$ introduced in \cite{Pan09-ord}. In Section 3 we present a concept of uniform distribution in infinite-dimensional rectangles for calculation of Riemann integrals for continuous functions over such rectangles(cf. \cite{Pan11-2}). In Section 4 we present Change of Variables Formula for $\alpha$-ordinary Lebesgue measure in $R^{\infty}$ established in \cite{Pan09}. In Section 5 we give a representation of the Dirac delta function in $ \mathcal{C}(R^{\infty})$ in terms of infinite-dimensional Lebesgue measures and consider  some properties of this functional.

\section{On ordinary and standard  Lebesgue measures  in $R^{\infty}$}

The problem of the existence of an  analog of the
Lebesgue measure for the  vector space of all real-valued
sequences $R^{\infty}=\prod_{i=1}^{\infty}R$
equipped with Tychonoff topology was discussed in \cite{Pan09-ord}.

R. Baker \cite{Bak91}~  firstly introduced the notion of ``Lebesgue
measure" in $R^{\infty}$ as follows:  a measure $\lambda$
being the completion of a translation invariant Borel measure
   in $R^{\infty}$  is called a ``Lebesgue measure" in  $R^{\infty}$ if
  for any measurable rectangle $\prod_{i=1}^{\infty}(a_i,b_i),~ -\infty < a_i < b_i <+\infty $ with $0 \le
\prod_{i=1}^{\infty}(b_i-a_i) < + \infty$, the following equality
$$\lambda\Big(\prod_{i=1}^{\infty}(a_i,b_i)\Big)=
\prod_{i=1}^{\infty}(b_i-a_i)$$ holds, where
$$\prod_{i=1}^{\infty}(b_i-a_i):=\lim_{n \to
\infty}\prod_{i=1}^{n}(b_i-a_i).$$

Subsequently, R. Baker \cite{Bak04} extended  his notion of
``Lebesgue measure" in $R^{\infty}$ as follows : a
measure $\lambda$ being the completion of a translation invariant
Borel measure
   on $R^{\infty}$  is called a ``Lebesgue measure"  if
  for any measurable rectangle $\prod_{i=1}^{\infty}R_i,~ R_i
\in \mathcal{B}(R)$~ with $0 \le
\prod_{i=1}^{\infty}m(R_i) < \infty$, the following equality
$$\lambda\Big(\prod_{i=1}^{\infty}R_i\Big)= \prod_{i=1}^{\infty}m(R_i)$$
holds, where $m$ denotes  a linear  Lebesgue measure in
$R$.

To propose a new concept of  ``Lebesgue measure" in
$R^{\infty}$, in \cite{Pan09-ord} main  attention  has been  attracted  to  the following
two simple facts:

\begin{fact}
Let $\mu$ be  a probability measure defined on a measure space
$(E,S)$. Then the product measure $\mu^{N}$ defined on
$(E^{N},S^{N})$ has the following essential
property: if $f$ is any permutation of $N$ and
$A_f((x_k)_{k \in N})=(x_{f(k)})_{k \in N}$ for
$(x_k)_{k \in N} \in  E^{N}$, then
$\mu^{N}(A_f(X))=\mu^{N}(X)$ for every $X \in
S^{N}$.
\end{fact}

\begin{fact}
The $n$-dimensional Lebesgue measure $\ell_n$ in $R^n$
has the following property: if $f$ is any permutation of $\{1,
\dots, n\}$ and
$$A_f((x_k)_{1 \le k \le n})=(x_{f(k)})_{1 \le k \le n}~((x_k)_{1 \le k \le n} \in R^n),$$ then
$\ell_n(A_f(X))=\ell_n(X)$ for every $X \in
\mathcal{B}(R^n).$
\end{fact}

In view of these facts  one can say that Baker's  measures
\cite{Bak91}, \cite{Bak04} have no  essential property of a
product - measure to be an invariant under the group  of all
canonical permutations \footnote{Let $f$ be any permutation of
$N$. A mapping $A_f:  {R^{\infty}}\to {R^{\infty}}$ defined by
$A_f((x_k)_{k \in N})=(x_{f(k)})_{k \in N}$ for
$(x_k)_{k \in N} \in R^{\infty}$ is called a
canonical permutation of $R^{\infty}.$} of
$R^{\infty}$.

Indeed, if we consider the following infinite-dimensional
rectangular set
$$X=\prod_{k=1}^{\infty}[0, e^{\frac{(-1)^k}{k}}],$$
then for every non-zero real number $a$ there exists a permutation
$f_{a}$ of  $N$ such that $\lambda(A_{f_a}(X))=a$, where
$\lambda$ is any Baker's  measure \cite{Bak91}, \cite{Bak04}.

To introduce  new concepts of the Lebesgue measure in
$R^{\infty}$, the following  definitions were introduced  in  \cite{Pan09-ord}:

\medskip

\begin{defn} Let $(\beta_j)_{j \in N} \in [0, +\infty]^N$.
We say that a number $\beta \in [0, +\infty]$ is an ordinary
product of numbers $(\beta_j)_{j \in N}$ if
$$
\beta=\lim_{n \to \infty}\prod_{i=1}^n\beta_i.
$$
An ordinary product of numbers $(\beta_j)_{j \in N}$ is
denoted by ${\bf (O)}\prod_{i \in N}\beta_i$.
\end{defn}

\begin{defn} Let $(\beta_j)_{j \in N} \in [0, +\infty]^N$.
A standard  product of the family of numbers $(\beta_i)_{i \in
N}$ is denoted by ${\bf (S)}\prod_{i \in
N}\beta_i$ and  defined as follows:

~${\bf (S)}\prod_{i \in N}\beta_i=0$ if ~$ \sum_{i \in
{N^{-}}}\ln(\beta_i)=-\infty$, where
$N^{-}=\{i:ln(\beta_i)<0\}$ ~\footnote{We set $ln(0)=-\infty$},
and ${\bf (S)}\prod_{i \in N}\beta_i=e^{\sum_{i \in
N}\ln(\beta_i)}$ if $\sum_{i \in
{N^{-}}}\ln(\beta_i) \neq -\infty$.
\end{defn}

Let $\alpha=(n_k)_{k \in N} \in (N \setminus
\{0\})^{N}$. We set
$$
F_0=[0,n_0] \cap N ,~F_1=[n_0\!+\!1,n_0\!+\!n_1]\cap
N,~\dots, F_k=[n_0\!+\cdots+ n_{k-1}\!+\!1,n_0\!+\cdots
+\!n_k]\cap N, \dots\,.
$$

\begin{defn}
We say that a number $\beta \in [0, +\infty]$ is an ordinary
$\alpha$-product of numbers $(\beta_i)_{i \in N}$ if
$\beta$ is   an ordinary product of numbers $(\prod_{i  \in
F_k}\beta_i)_{k \in N}$. An ordinary $\alpha$-product of
numbers $(\beta_i)_{i \in N}$ is denoted by ${\bf
(O,\alpha)}\prod_{i \in N}\beta_i$.
\end{defn}

\begin{defn}
We say that a number $\beta \in [0, +\infty]$ is a standard
$\alpha$-product of numbers $(\beta_i)_{i \in N}$ if
$\beta$  is a standard product of numbers $(\prod_{i  \in
F_k}\beta_i)_{k \in N}$. A standard $\alpha$-product of
numbers $(\beta_i)_{i \in N}$ is denoted by ${\bf
(S,\alpha)}\prod_{i \in N}\beta_i$.
\end{defn}

\begin{defn}
Let $\alpha=(n_k)_{k \in N} \in (N \setminus
\{0\})^{N}$. Let $\mathcal{(\alpha)OR} $ be the class of
all infinite-dimensional measurable $\alpha$-rectangles  $R =
\prod_{i \in N}R_i (R_i \in
\mathcal{B}(R^{n_i})$) for which an ordinary product of
numbers $(m^{n_i}(R_i))_{i \in N}$ exists and is finite.
We say that a  measure  $\lambda$ being the completion of a
translation-invariant Borel measure  is an ordinary
$\alpha$-Lebesgue measure in ${R^{\infty}}$(or, shortly,
O$(\alpha)$LM)  if ~for every $R  \in \mathcal{(\alpha)OR} $ we
have
$$
\lambda(R)={\bf (O)}\prod_{k \in N}m^{n_k}(R_k).
$$
\end{defn}

\begin{defn}
Let $\alpha=(n_k)_{k \in N} \in (N \setminus
\{0\})^{N}$. Let $\mathcal{(\alpha)SR} $ be the class of
all infinite-dimensional measurable $\alpha$-rectangles  $R =
\prod_{i \in N}R_i (R_i \in
\mathcal{B}(R^{n_i}))$ for which a standard product of
numbers $(m^{n_i}(R_i))_{i \in N}$ exists and is finite.
We say that a  measure  $\lambda$ being the completion of a
translation-invariant Borel measure  is a standard
$\alpha$-Lebesgue measure in ${R^{\infty}}$(or, shortly,
S$(\alpha)$LM)  if ~for every $R  \in \mathcal{(\alpha)SR} $ we
have
$$
\lambda(R)={\bf (S)}\prod_{k \in N}m^{n_k}(R_k).
$$
\end{defn}

\begin{proposition}(\cite{Pan09-ord}, Proposition 1, p. 212)
Note that for every $\alpha=(n_k)_{k \in N} \in
(N \setminus \{0\})^{N}$ the following strict
inclusion
$$
 \mathcal{(\alpha)OR} \subset \mathcal{(\alpha)SR}$$ holds.
\end{proposition}

The presented approach  gives us a possibility to
construct  such  translation-invariant  Borel measures in $
R^{\infty}$ which are different from the Baker measures
\cite{Bak04} in the sense that it does not apply the metric
properties of $R^{\infty}$. It is an adaptation of
a construction from general measure theory which allows us to
construct interesting examples of analogs of a Lebesgue measure on
the entire space.

\medskip

Let $(E,S)$ be a measurable space and let $\mathcal{R}$ be any
subclass of the $\sigma$-algebra $S$. Let $(\mu_B)_{B \in
\mathcal{R}}$ be such a family of $\sigma$-finite measures that
 for $B \in \mathcal{R}$ we have $\mbox{dom}(\mu_B)=S \cap \mathcal{P}(B)$,
 where $\mathcal{P}(B)$ denotes the power set of the set $B$.

\begin{defn}
A family  $(\mu_B)_{B \in \mathcal{R}}$ is called to be consistent
if
$$
(\forall X)(\forall B_1,B_2)(X \in S~\&~B_1,B_2 \in \mathcal{R}
\rightarrow \mu_{B_1}(X \cap B_1 \cap B_2)=\mu_{B_2}(X \cap B_1
\cap B_2)).
$$
\end{defn}

The following assertion plays a key role for construction of new translation-invariant measures.

\begin{lem} (\cite{Pan09-ord}, Lemma 1, p. 213)
Let $(\mu_B)_{B \in \mathcal{R}}$ be a consistent family of
$\sigma$-finite measures. Then there exists a measure
$\mu_{\mathcal{R}}$ on $(E,S)$ such that

{\rm (i)}~$\mu_{\mathcal{R}}(B)=\mu_B(B)$ for every $B \in
\mathcal{R}$;

{\rm (ii)}~if there exists a non-countable family of pairwise
disjoint sets $\{B_i~:~i \in I\} \subseteq \mathcal{R}$ such that
$0<\mu_{B_i}(B_i)<\infty$, then the measure $\mu_{\mathcal{R}}$ is
non-$\sigma$-finite;

{\rm (iii)}~if  $G$ is a group of measurable transformations of
$E$ such that $G(\mathcal{R})=\mathcal{R}$ and
$$
(\forall B)(\forall X)(\forall g)\big(\big(B \in \mathcal{R}~\& X
\in S \cap \mathcal{P}(B)~\&~ g \in G \big) \rightarrow
\mu_{g(B)}(g(X))=\mu_{B}(X)\big),
$$
where $\mathcal{P}(B)$ denotes a power set of the set $B$, then
the measure $\mu_{\mathcal{R}}$ is $G$-invariant.
\end{lem}

\begin{lem}(\cite{Pan09-ord}, Lemma 2. p. 216)
Let $\alpha=(n_i)_{i \in N} \in (N \setminus
\{0\} )^{N}$. We set $\mathcal{R}=\mathcal{(\alpha)OR}$.
Suppose that $R=\prod_{i \in N}R_i \in \mathcal{R}$ for
which $R_i \in \mathcal{B}(R^{n_i})$ for $i \in
N$.

For  $X \in \mathcal{B}(R)$, we set $\mu_R(X)=0$ if $$ {\bf
(O)}\prod_{i \in N}m^{n_i}(R_i)=0,$$ and
$$
\mu_R(X)={\bf (O)}\prod_{i \in N}m^{n_i}(R_i) \times
\Big(\prod_{i \in
N}\frac{{m^{n_i}}_{R_i}}{m^{n_i}(R_i)}\Big)(X)
$$
otherwise, where $\frac{{m^{n_i}}_{R_i}}{m^{n_i}(R_i)}$ is a Borel
probability measure defined on $R_i$ as follows
$$
(\forall X)\Big(X \in \mathcal{B}(R_i) \rightarrow
\frac{{m^{n_i}}_{R_i}}{m^{n_i}(R_i)}(X)=\frac{m^{n_i}(Y \cap
R_i)}{m^{n_i}(R_i)}\Big).
$$

 Then the  family of measures $(\mu_{R})_{R \in
\mathcal{R}}$ is consistent.
\end{lem}

\begin{lem}(\cite{Pan09-ord}, Lemma 2. p. 217)
Let $\alpha=(n_i)_{i \in N} \in (N \setminus
\{0\})^{N}$. We set $\mathcal{R}=\mathcal{(\alpha)SR}$.
Suppose that $R=\prod_{i \in N}R_i \in \mathcal{R}$ for
which $R_i \in \mathcal{B}(R^{n_i})$ for $i \in
N$ and $R \in \mathcal{(\alpha)SR}$.

For  $X \in \mathcal{B}(R)$, we set $\mu_R(X)=0$ if $$ {\bf
(S)}\prod_{i \in N}m^{n_i}(R_i)=0,$$ and
$$
\mu_R(X)={\bf (S)}\prod_{i \in N}m^{n_i}(R_i) \times
\Big(\prod_{i \in
N}\frac{{m^{n_i}}_{R_i}}{m^{n_i}(R_i)}\Big)(X)
$$
otherwise, where $\frac{{m^{n_i}}_{R_i}}{m^{n_i}(R_i)}$ is a Borel
probability measure defined in $R_i$ as follows
$$
(\forall X)\Big(X \in \mathcal{B}(R_i) \rightarrow
\frac{{m^{n_i}}_{R_i}}{m^{n_i}(R_i)}(X)=\frac{m^{n_i}(Y \cap
R_i)}{m^{n_i}(R_i)}\Big).
$$

 Then the  family of measures $(\mu_{R})_{R \in
\mathcal{R}}$ is consistent.
\end{lem}

\medskip

Next two theorems are corollaries of  Lemmas 2.12--2.13.

\begin{thm}(\cite{Pan09-ord}, Theorem 1. p. 217)
For every  $\alpha=(n_i)_{i \in N} \in (N
\setminus \{0\})^{N}$, there exists a Borel measure
$\mu_{\alpha}$ in $R^{\infty}$ which is O$(\alpha)$LM.
\end{thm}

\begin{thm}(\cite{Pan09-ord}, Theorem 1. p. 218)
For every  $\alpha=(n_i)_{i \in N} \in (N
\setminus \{0\})^{N}$, there exists a Borel measure
$\nu_{\alpha}$ in $R^{\infty}$ which is S$(\alpha)$LM.
\end{thm}

Let $\mu_1$ and $\mu_2$ be two measures defined on the measurable
space $(~{E},~{S})$.

\begin{defn}[\cite{Hal50}, p.~124]
We say that the $\mu_1$ is absolutely continuous with respect to
the $\mu_2$, in symbols $\mu_1 \ll \mu_2$, if
$$
(\forall X)(X \in ~{S}~\&~\mu_2(X)=0 \rightarrow
\mu_1(X)=0).
$$
\end{defn}

\begin{defn}[\cite{Hal50}, p.~126]
Two measures $\mu_1$ and $\mu_2$ for which both  $\mu_1 \ll \mu_2$
and $\mu_2 \ll \mu_1$ are called equivalent, in symbols $\mu_1
\equiv \mu_2$.
\end{defn}

We have the following  assertion.

\begin{thm} (\cite{Pan09-ord}, Theorem 3. p. 217)
For every  $\alpha=(n_i)_{i \in N} \in (N
\setminus \{0\})^{N}$, we have $\nu_{\alpha} \ll
\mu_{\alpha}$ and the measures $\nu_{\alpha}$ and $\mu_{\alpha}$
are not equivalent.
\end{thm}

\begin{rem}
Note that  the $\mu_{\alpha}$ coincides with Baker's measure
\cite{Bak04} for $\alpha=(1,1, \dots)$. By Lemmas 2.12 and 2.13 we
can get the construction of  Baker's  measure \cite{Bak91}. In
this direction we must consider a class $\mathcal{R}_B$ of all
measurable rectangles $\prod_{i=1}^{\infty}(a_i,b_i),~ -\infty <
a_i < b_i <+\infty $ for which $0 \le ({\bf O})\prod_{i \in
N}(b_i-a_i) < + \infty$. Since $\mathcal{R}_B$ is
translation-invariant and the family of measures  $(\mu_R)_{R \in
\mathcal{R}_B}$ is consistent as a subfamily of the consistent
family of measures constructed in Lemma 2.12, we claim that Baker's
measure \cite{Bak91} coincides with the measure
$\lambda_{\mathcal{R}_B}$.
\end{rem}

\begin{defn}
Let $\alpha=(n_i)_{i \in N} \in (N \setminus
\{0\})^{N}$ such that $n_i=n_j$ for every $i,j \in
N$. We set $F_i=(a_1^{(i)},\dots,a_{n_0}^{(i)})$ for
every $i \in N$(see, notations introduced  before
Definition 2.5. Let $f$ be any permutation of $N$ such
that for every $i \in N$ there exists $j \in N$
such that $f(a_k^{(i)})=a_k^{(j)}$ for
 $1 \le k \le n_0$. Then a map $A_f:
R^{\infty}\to R^{\infty}$ defined by
$A_f((z_k)_{k \in N})=(z_{f(k)})_{k \in N}$ for
$(z_k)_{k \in N}\in R^{\infty}$, is called a
canonical $\alpha$-permutations of $R^{\infty}$.
\end{defn}

A group of transformations generated by all $\alpha$-permutations
and shifts of  $R^{\infty}$, is denoted by
$\mathcal{G}_{\alpha}$.

\begin{cor}
For every  $\alpha=(n_i)_{i \in N} \in (N
\setminus \{0\})^{N}$ for which  $n_i=n_j ( i,j \in
N)$,    the measure $\nu_{\alpha}$ is
$\mathcal{G}_{\alpha}$-invariant.
\end{cor}

One can easily get the validity of the following propositions.

\begin{proposition}  (\cite{Pan09-ord}, Proposition 2, p. 219)
For every  $\alpha=(n_i)_{i \in N} \in (N
\setminus \{0\})^{N}$ there exists $\beta \in (N
\setminus \{0\})^{N}$ such that $\mu_{\alpha}$ and
$\mu_{\beta}$ are different.
\end{proposition}

\begin{proposition}(\cite{Pan09-ord}, Proposition 2, p. 220)
For every  $\alpha=(n_i)_{i \in N} \in (N
\setminus \{0\})^{N}$ there exists $\beta \in (N
\setminus \{0\})^{N}$ such that $\nu_{\alpha}$ and
$\nu_{\beta}$ are different.
\end{proposition}

\section{On uniformly  distributed  sequences  of increasing family of finite sets in infinite-dimensional rectangles}

\vspace{.08in} \noindent Let $s_1, s_2, s_3, \dots$ be a uniformly
distributed in an interval $[a, b]$ (see, for example \cite{KuiNie74}. Setting $Y_n=\{s_1, s_2, s_3,
\dots , s_n\}$ for $n \in N$, the $(Y_n)_{n \in N}$ will be such
an increasing sequence of finite subsets of the $[a, b]$ that, for
any subinterval $[c, d]$ of the $[a, b]$, the following equality
$$
\lim_{n \to \infty}\frac{\#(Y_n \cap
[c,d])}{\#(Y_n)}=\frac{d-c}{b-a}
$$
will be valid.

This remark raises the following


\begin{defn}

An increasing sequence $(Y_n)_{n \in N}$ of  finite subsets of the
$[a, b]$ is said to be equidistributed or uniformly distributed in
an interval $[a, b]$ if, for any subinterval $[c, d]$ of the $[a,
b]$, we have
$$
\lim_{n \to \infty}\frac{\#(Y_n \cap
[c,d])}{\#(Y_n)}=\frac{d-c}{b-a}.
$$
\end{defn}

\begin{defn}
  Let $\prod_{k \in N}[a_k,b_k] \in
\mathcal{R}$. A set $U$ is called an elementary rectangle in the
$\prod_{k \in N}[a_k,b_k]$  if it admits the following
representation
 $$
 U=\prod_{ k=1}^m ][c_k, d_k][ \times \prod_{ k \in N \setminus \{1,\dots,
 m\}}[a_k,b_k],
 $$
where   $a_k \le c_k < d_k \le b_k$ for $1 \le k \le m.$
\end{defn}

It is obvious that
$$ \lambda(U)=\prod_{
k=1}^m(d_k-c_k)\times \prod_{ k=m+1}^{\infty}(b_k-a_k),
$$
for the elementary rectangle $U$.

\begin{defn}
 An increasing
sequence $(Y_n)_{n \in N}$ of finite subsets of the
infinite-dimensional  rectangle $\prod_{k \in N}[a_k,b_k] \in
\mathcal{R}$ is said to be uniformly distributed in the $\prod_{k
\in N}[a_k,b_k]$
 if for every elementary  rectangle $U$ in the $\prod_{k \in N}[a_k,b_k[$
 we have
$$
\lim_{n \to \infty}\frac{\#(Y_n \cap U)}{\#(Y_n) }=
\frac{\lambda(U)}{\lambda(\prod_{k \in N}[a_k,b_k[)}.
$$
\end{defn}


\begin{thm} (\cite{Pan11-2}, Theorem 3.1, p.328)
 Let $\prod_{k \in N}[a_k,b_k]
\in \mathcal{R}$.
 Let $(x_n^{(k)})_{n \in N}$ be uniformly distributed
in the interval $[a_k,b_k]$ for $ k \in N$. We set
$$Y_n=\prod_{k=1}^n (\cup_{j=1}^n x_j^{(k)})\times \prod_{k \in N \setminus \{1,\dots,
 n\}}\{x_0^{(k)}\}.
 $$
 Then $(Y_n)_{n \in N}$ is uniformly
distributed in the rectangle $\prod_{k \in N}[a_k,b_k]$.
\end{thm}

\begin{defn}
 Let $\prod_{k \in N}[a_k,b_k] \in
\mathcal{R}$. A family of pairwise disjoint elementary rectangles
$\tau=(U_k)_{1 \le k \le n}$ of the $\prod_{k \in N}[a_k,b_k]$ is
called Riemann partition of the $\prod_{k \in N}[a_k,b_k]$ if
$\cup_{1 \le k \le n}U_k=\prod_{k \in N}[a_k,b_k]$.
\end{defn}

\begin{defn}

Let $\tau=(U_k)_{1 \le k \le n}$ be Riemann partition of the
$\prod_{k \in N}[a_k,b_k]$.  Let $\ell(Pr_i(U_k))$ be a length of
the $i$-th projection $Pr_i(U_k)$ of the $U_k$ for $i \in N$. We
set
$$
d(U_k)=\sum_{i \in
N}\frac{\ell(Pr_i(U_k))}{2^i(1+\ell(Pr_i(U_k)))}.
$$
It is obvious that $d(U_k)$ is a diameter of the elementary
rectangle $U_k$ for $k \in N$ with respect to Tikhonov metric
$\rho$ defined as follows
$$ \rho((x_k)_{k \in N},
(y_k)_{k \in N})=\sum_{k \in N}\frac{|x_k-y_k|}{2^k(1+|x_k-y_k|)}
$$ for $(x_k)_{k \in N}, (y_k)_{k \in N} \in {\bf
R}^{\infty}.$

A number $d(\tau)$, defined by
$$d(\tau)=\max \{ d(U_k) : 1 \le k \le
n \}
$$
is called mesh or norm  of the Riemann partition $\tau$.
\end{defn}


\begin{defn}

 Let
$\tau_1=(U_i^{(1)})_{1 \le i \le n}$ and $\tau_2=(U_j^{(2)})_{1
\le j \le m}$ be Riemann partitions of the $\prod_{k \in
N}[a_k,b_k]$.  We say that $\tau_2 \le \tau_1$ iff
$$
(\forall j)((1 \le j \le m) \rightarrow (\exists i_0)(1 \le i_0
\le n ~\&~ U_j^{(2)} \subseteq U_{i_0}^{(1)})).
$$
\end{defn}

\begin{defn}

 Let $f$ be a real-valued bounded function
defined on the $\prod_{i \in N}[a_i,b_i]$. Let $\tau=(U_k)_{1 \le
k \le n}$ be Riemann partition of the $\prod_{k \in N}[a_k,b_k]$
and $(t_k)_{1 \le k \le n}$ be a sample such that, for each $k$,
$t_k \in U_k$. Then

(i) a sum  $\sum_{k=1}^n f(t_k) \lambda(U_k) $ is called Riemann
sum of the $f$ with respect to Riemann  partition $\tau=(U_k)_{1
\le k \le n}$  together with sample $(t_k)_{1 \le k \le n}$;

(ii) ~ a sum $S_{\tau}=\sum_{k=1}^n M_k \lambda(U_k)$ is called
the upper Darboux sum with respect to Riemann  partition $\tau$,
where $M_k=\sup_{x \in U_k}f(x)(1 \le k \le n)$;

(ii) ~ a sum $s_{\tau}=\sum_{k=1}^n m_k \lambda(U_k)$ is called
the lower Darboux sum with respect to Riemann  partition $\tau$,
where $m_k=\inf_{x \in U_k}f(x)(1 \le k \le n)$.
\end{defn}


\begin{defn}

 Let $f$ be a
real-valued bounded function defined on $\prod_{i \in
N}[a_i,b_i[$. We say that the $f$ is Riemann-integrable  on
$\prod_{i \in N}[a_i,b_i]$ if there exists
 a real number $s$ such that for every positive real number
 $\epsilon$ there exists a real number $ \delta >0$ such that,
 for every Riemann partition $(U_k)_{1 \le k \le n}$ of the $\prod_{k \in N}[a_k,b_k]$
 with $d(\tau)<\delta$ and for every sample $(t_k)_{1 \le k
\le n}$,  we have
$$
\big|\sum_{k=1}^n f(t_k) \lambda(U_k)-s\big|<\epsilon.
$$
The number $s$ is called Riemann integral and is denoted by
$$(R)\int_{\prod_{k \in N}[a_k,b_k]}f(x)d\lambda(x).$$
\end{defn}


\begin{defn}

A function $f$ is called a step function on $\prod_{k \in
N}[a_k,b_k]$ if it can be written~as
$$
f(x)=\sum_{k=1}^nc_k\mathcal{X}_{U_k}(x),
$$
where $\tau=(U_k)_{1 \le k \le n}$ is any Riemann partition of the
$\prod_{k \in N}[a_k,b_k]$, $c_k \in R$ for $1 \le k \le n$ and
$\mathcal{X}_{A}$ is the indicator function of the $A$
\end{defn}

\begin{thm} (\cite{Pan11-2}, Theorem 3.2, p.331)
Let $f$ be a continuous function on $\prod_{k \in N}[a_k,b_k]$
with respect to Tikhonov metric $\rho$. Then the  $f$  is
Riemann-integrable on $\prod_{k \in N}[a_k,b_k]$.
\end{thm}

Let denote by $\mathcal{C}(\prod_{k \in N}[a_k,b_k])$
a class of all continuous (with respect to Tikhonov topology)
real-valued functions on $\prod_{k \in N}[a_k,b_k]$.


\begin{thm} (\cite{Pan11-2}, Theorem 3.4, p.336)
 For $\prod_{i \in
N}[a_i,b_i] \in \mathcal{R}$, let  $(Y_n)_{n \in N}$ be an
increasing family its finite subsets. Then $(Y_n)_{n \in N}$ is
uniformly distributed in the $\prod_{k \in N}[a_k,b_k]$ if and
only if for every  $f \in \mathcal{C}(\prod_{k \in N}[a_k,b_k])$
the following equality
$$
\lim_{n \to \infty}\frac{\sum_{y \in
Y_n}f(y)}{\#(Y_n)}=\frac{(R)\int_{\prod_{k \in
N}[a_k,b_k]}f(x)d\lambda(x)}{\lambda(\prod_{i \in
N}[a_i,b_i])}
$$ holds.
\end{thm}

\section{Change of variable formula for the $\alpha$-ordinary  Lebesgue measure in  $~{R^N}$}

Let  $R^n(n>1)$ be an $n$-dimensional Euclidean space and
let $\mu_n$ an $n$-dimensional standard Lebesgue measure on
$R^n$. Further, let $T$ be a linear $\mu_n$-measurable
transformation of $R^n$.

 It is obvious that $\mu_n T^{-1}$ is absolutely continuous
with respect to $\mu_n$, and ~ there exists a non-negative
$\mu_n$-measurable function $\Phi$ on $R^n$ such that
$$
\mu_n (T^{-1}(X))=\int_X \Phi(y)d \mu_n(y)
$$
for every $\mu_n$-measurable subset $X$ of $R^n$.

  The function $\Phi$ plays the role of the Jacobian $J(T^{-1})$
of the transformation $T^{-1}$(or, rather the absolute value of
the Jacobian)(see, e.g., \cite{Hal50}) in the theory of
transformations of multiple integrals. It is clear that
$J(T^{-1})$ coincides with a Radon-Nikodym derivative $\frac{d
\mu_n T^{-1}}{d\mu_n}$, which is unique a.e. with respect to
$\mu_n$.

 It is clear that
$$
\frac{d\mu_nT^{-1}}{d\mu_n}(x)=\lim_{k \to
\infty}\frac{\mu_n(T^{-1}(U_k{(x)}))}{\mu_n(U_k(x))}(\mu_n-
\mbox{a.e.}),
$$
where $U_k(x)$ is a spherical neighborhood with the center in $x
\in R^n $ and radius ${\bf r}_k >0$ so that $\lim_{k \to
\infty} r_k=0$. The class of such spherical neighborhoods generate
so-called Vitali differentiability class of subsets which allows
us to calculate the Jacobian $J(T^{-1})$ of the transformation
$T^{-1}$.

If we consider a vector space of all real-valued sequences
$~{R^N}$(equipped with Tychonoff topology), then we observe
that  for the infinite-dimensional Lebesgue measure \cite{Bak91}
(or \cite{Bak04}) defined in $~{R^N}$ there does not exist
any Vitali system of differentiability,  but in spite of
non-existence of such a system the inner structure of this measure
allows us to define a form of the Radon-Nikodym derivative defined
by any linear transformation of $~{R^N}$. In order to show
it, let consider the following

\begin{example}
Let $ \mathcal{R}_1$ be the class of all infinite dimensional
rectangles $R  \in B(~{R^N})$ of the form

$$R = \prod_{i=1}^{\infty}R_i, ~~ R_i=(a_i,b_i), -\infty < a_i \le  b_i < +\infty,$$
such that $$0 \le \prod_{i=1}^{\infty}(b_i-a_i)<\infty.$$ Let
$\tau_1$ be the set function on $\mathcal{R}_1$ defined by
$$
\tau_1(R)=\prod_{i=1}^{\infty}(b_i-a_i).
$$

 R.~Baker \cite{Bak91}  proved that the functional $\lambda_1$ defined by
$$
( \forall X)\Big(X \in  B(~{R^N}) \rightarrow \lambda_1(X) =
\inf \Big\{ \sum_{j=1}^{\infty} \tau_1(R_j) : R_j \in
\mathcal{R}_1~ \&~ X \subseteq  \cup_{j=1}^{\infty}R_j \Big\}\Big)
$$
is a quasi-finite translation-invariant Borel measure in
$~{R^N}$.

The following change of variable formula has been established in
 \cite{Bak91} (cf. p.~1029): {\it  Let $T^n : R^n \to R^n,~ n >
1$, be a linear transformation with Jacobian $\Delta \neq0,$ and
let $T^{N} : ~{R^{N}} \to ~{R^{N}}$ be the
map defined by $$T^{N}(x) = (T^n(x_1, \dots, x_n),
x_{n+1}, x_{n+2}, \dots), ~x = (x_i)_{i \in N} \in ~{R^N}.
$$ Then for each $ E \in
\mathcal{B}(~{R^{N}}),$ ~we have
$$\lambda_1(T^{N}(E)) =
|\Delta|\lambda_1(E).$$}
\end{example}

\begin{thm}
Let $\alpha=(n_i)_{i \in N}$ be the sequence of non-zero natural
numbers  and $\mu_{\alpha}$ is $\mathcal{O}(\alpha)LM$. Further, let $T^{n_i} : R^{n_i} \to R^{n_i},
i \ge 1$, be a  family of linear transformation with Jacobians
$\Delta_i \neq 0$ and $0 < \prod_{i=1}^{\infty}\Delta_i < \infty$.
Let $T^{N} : ~{R^{N}} \to ~{R^{N}}$ be the map defined
by
$$T^{N}(x) =
(T^{n_1}(x_1, \dots, x_{n_1}),T^{n_2}(x_{n_1+1}, \dots,
x_{n_1+n_2}), \dots),$$ where $x = (x_i)_{i \in N} \in
~{R^{N}}$. Then for each $ E \in
\mathcal{B}(~{R^{N}}),$ ~we have
$$\mu_{\alpha}(T^{N}(E)) =
\Big(\prod_{i=1}^{\infty}\Delta_i\Big) \mu_{\alpha}(E).$$
\end{thm}

\begin{rem}
Theorem 4.5 is  change of variable formula for the
$\alpha$-ordinary  Lebesgue measure.  It extends change of
variable formula for Baker's measure considered in Example 4.1.
Indeed, let $T^n : R^n \to R^n, n
> 1$, be a linear transformation with Jacobian $\Delta \neq 0$.
Let $n_1=n$ and $n_i=1$ for $i > 1$, that is $\alpha=(n,1,1,\cdots)$. Further, we set $T^{n_1}=T^n$
and $T^{n_k}=I$, where $I : R \to R$ is an
identity transformation of $R$ defined by $I(x)=x$ for $x
\in R$.
\end{rem}

Let a map $T^{N} : ~{R^{N}} \to ~{R^{N}}$ be
defined by
$$T^{N}(x) = (T^n(x_1, \dots, x_n), x_{n+1}, x_{n+2}, \dots),~ x
= (x_i)_{i \in N} \in  ~{R^N} .$$

Then, by Theorem 4.5,  for $T^{N}$ and for each $ E
\in \mathcal{B}(~{R^{N}}),$ ~we have
$$\lambda_1(T^{N}(E)) =\mu_{\alpha}(T^{N}(E))=|\Delta|\mu_{\alpha}(E)=
|\Delta|\lambda_1(E).$$

\section{Concept of the  Dirac delta function in $\mathcal{C}(R^{\infty}$)}

\begin{lem}(Intermediate value theorem) Let $f$ be a continuous
function on $\prod_{k \in N}[a_k,b_k]$. Suppose that $\max\{f(x):x
\in \prod_{k \in N}[a_k,b_k]\}=M$ and $\min\{f(x):x \in \prod_{k
\in N}[a_k,b_k]\}=m$. Let $u \in [m,M]$. Then there is $c \in
\prod_{k \in N}[a_k,b_k]$ such that $f(c)=u$.
\end{lem}

{\bf Proof.}  Let $(y_k)_{k \in N} \in \prod_{k \in N}[a_k,b_k]$ be such
sequence that $f((y_k)_{k \in N})=M$.

Let $Z^*=(z_k)_{k \in N} \in \prod_{k \in N}[a_k,b_k]$ be such
a sequence that $f((z_k)_{k \in N})=m$.

Let consider a function $g(t)=f((z_k)_{k \in N}+t((y_k)_{k \in
N}-(z_k)_{k \in N}))$ on $[0,1]$. This function is well defined
because $(z_k)_{k \in N}+t((y_k)_{k \in N}-(z_k)_{k \in N})$ is in
$\prod_{k \in N}[a_k,b_k]$ for each $t \in [0,1]$. It is obvious
that $g_{max}=g(1)=M$ and $g_{min}=g(0)=m$. Using Intermediate
value theorem for a real valued function $g$ on $[0,1]$ there is
$t_0 \in [0,1]$ such that $g(t_0)=u$. Setting $c=(z_k)_{k \in
N}+t_0((y_k)_{k \in N}-(z_k)_{k \in N})$, we end the proof of the
lemma.

$$ \eqno \diamondsuit$$

Let $\lambda$ be Baker  measure  in $R^{\infty}$. For $\epsilon>0$, we set $$a_k(\epsilon)=\frac{e^{-\frac{1}{2^k \epsilon}}}{2}$$
and
$$\Delta_{\epsilon}=\prod_{k=1}^{\infty}[-a_k(\epsilon),a_k(\epsilon)].$$

Note that the diameter of the set $\Delta_{\epsilon}$ is
calculated by
$$
\mbox{diam}(\Delta_{\epsilon})=\sum_{i \in
N}\frac{|2a_k(\epsilon))|}{2^i(1+|2a_k(\epsilon)|)}=\sum_{i \in
N}\frac{e^{-\frac{1}{2^k \epsilon}}}{2^i(1+e^{-\frac{1}{2^k
\epsilon}})}.
$$

\begin{lem} $\lim_{\epsilon \to 0+}\mbox{diam}(\Delta_{\epsilon})=0$.
\end{lem}

{\bf Proof.}

For $\sigma >0$ there is $n_{\sigma}\in N$ such that
$$\sum_{i=n_{\sigma}}^{\infty}2^{-i}< \frac{\sigma}{2}.
$$

Since $\lim_{\epsilon \to 0^{+}}\frac{e^{-\frac{1}{2^k
\epsilon}}}{2^i(1+e^{-\frac{1}{2^k \epsilon}})}=0$ for each $k \in
N$, we deduce that

$$\lim_{\epsilon \to  O+}\sum_{i=1}^{n_{\sigma}}\frac{e^{-\frac{1}{2^k
\epsilon}}}{2^i(1+e^{-\frac{1}{2^k \epsilon}})}=0.$$ The latter
relation means that  there is $\rho_{\sigma}>0$ such that
$$
\sum_{i=1}^{n_{\sigma}}\frac{e^{-\frac{1}{2^k
\epsilon}}}{2^i(1+e^{-\frac{1}{2^k \epsilon}})}<\frac{\sigma}{2}
$$
for all $\epsilon$  with $0<\epsilon< \rho_{\sigma}$.

Finally, for each $\sigma>0$,  $\rho_{\sigma}$ is such a positive
number that
$$
\mbox{diam}(\Delta_{\epsilon})=\sum_{i \in
N}\frac{|2a_k(\epsilon))|}{2^i(1+|2a_k(\epsilon)|)}=\sum_{i \in
N}\frac{e^{-\frac{1}{2^k \epsilon}}}{2^i(1+e^{-\frac{1}{2^k
\epsilon}}}\le $$
$$ \sum_{i=1}^{n_{\sigma}}\frac{e^{-\frac{1}{2^k
\epsilon}}}{2^i(1+e^{-\frac{1}{2^k
\epsilon}})}+\sum_{i=n_{\sigma}}^{\infty}2^{-i}\le
\frac{\sigma}{2}+\frac{\sigma}{2}=\sigma
$$
for each  $\epsilon $ with  with $0<\epsilon< \rho_{\sigma}$.

This ends the proof of the lemma.

For $y \in R^{\infty}$ we set
$\Delta_{\epsilon}(y)=\Delta_{\epsilon}+y.$

Since Tychonoff metric is translation invariant, by virtue of
Lemma 5.2  we deduce that
$$\lim_{\epsilon \to  O+}\mbox{diam}(\Delta_{\epsilon}(y))=0.$$

Note also that if $\epsilon^{(i)}>0$ for $i \in N$ and  $\lim_{i \to \infty} \epsilon^{(i)}=0$ then the equality
$$
\cap_{i \in N}\Delta_{\epsilon^{(i)}}(y)=\{y\}
$$
holds true for each $y \in R^{\infty}$.

\begin{lem}
Let $f$ be a continuous function on $R^{\infty}$. Then
the following formula
$$
\lim_{\epsilon \to  O+}\frac{1}{\lambda(\Delta_{\epsilon}(y))}
\int_{\Delta_{\epsilon}(y))}f(x) d\lambda(x)=f(y) \eqno(*)
$$
holds true for all $y \in R^{\infty}$.
\end{lem}

{\bf Proof.} If consider the restriction of $f$ on $\Delta_{\epsilon}(y)$
is also continuous. Let denote by $M_{\epsilon}(y)$ and
$m_{\epsilon}(y)$ maximum and minimum of the function of $f$ on
$\Delta_{\epsilon}(y)$. Hence we have
$$
m_{\epsilon}\times \lambda(\Delta_{\epsilon}(y)) \le
\int_{\Delta_{\epsilon}(y))}f(x) d\lambda(x)\le M_{\epsilon}\times
\lambda(\Delta_{\epsilon}(y))
$$
for each $\epsilon >0$. Equivalently, we have

$$
m_{\epsilon} \le  \frac{1}{\lambda(\Delta_{\epsilon}(y))}
\int_{\Delta_{\epsilon}(y))}f(x) d\lambda(x) \le M_{\epsilon}
$$
for each $\epsilon >0$.

By Lemma 5.1, there is $y_{\epsilon} \in
\Delta_{\epsilon}(y)$ such that
$$
\frac{1}{\lambda(\Delta_{\epsilon}(y))}
\int_{\Delta_{\epsilon}(y))}f(x) d\lambda(x) =f(y_{\epsilon}).$$
When one takes the limit when $\epsilon \to  O+$, then  $y_{\epsilon}$ tends to $y$,
and so
$$
\lim_{\epsilon \to  O+}\frac{1}{\lambda(\Delta_{\epsilon}(y))}
\int_{\Delta_{\epsilon}(y))}f(x) d\lambda(x)=\lim_{\epsilon \to O+}f(y_{\epsilon})=f(y).
$$

~

We have
$$
\lambda(\Delta_{\epsilon})=\prod_{k=1}^{\infty}(2a_k(\epsilon))=e^{-\sum_{k=1}^{\infty}\frac{1}{2^k \epsilon}}.
$$

We set $\eta_{\epsilon}(x)=e^{\sum_{k=1}^{\infty}\frac{1}{2^k \epsilon}}$ if $x \in
\Delta_{\epsilon}$ and $\eta_{\epsilon}(x)=0$, otherwise.

$\eta_{\epsilon}(x)$ is  called a nascent delta function.

The Dirac delta function $\delta(x)$, formally is defined by

$$
\delta(x)=\lim_{\epsilon \to  O+}\eta_{\epsilon}(x),
$$
which, of course, has no any reasonable sense.

 Let $f$ be a continuous real-valued function on
$R^{\infty}$. We define a Dirac delta integral as follows
$$(\delta)\int_{R^{\infty}}\delta(x) f(x)
d\lambda(x)= \lim_{\epsilon \to  O+}\int_{R^{\infty}}\eta_{\epsilon}(x) f(x) d\lambda(x).
$$

We define a Dirac delta functional $\delta
:~{C}(R^{\infty})\to R$ by
$$
\delta(f) =(\delta)\int_{R^{\infty}}\delta(x) f(x) d\lambda(x).
$$

The following assertion is valid.

\begin{thm} The Dirac delta functional $\delta$ is a linear functional such that $\delta(f)=f({\bf 0})$ for
each $f \in ~{C}(R^{\infty})$, where ${\bf 0}$
denotes the zero of $R^{\infty}$.
\end{thm}
{\bf Proof.}  We have
$$\delta(f)=(\delta)\int_{R^{\infty}}\delta(x) f(x)
d\lambda(x)= \lim_{\epsilon \to  O+}\int_{R^{\infty}}\eta_{\epsilon}(x) f(x) d\lambda(x)= $$
$$\lim_{\epsilon \to  O+}\int_{R^{\infty}}
\big[e^{\sum_{k=1}^{\infty}\frac{1}{2^k \epsilon}} \times \chi_{\Delta_{\epsilon}}(y)+0
\times \chi_{R^{\infty} \setminus \Delta_{\epsilon}}(y)\big]
 f(x) d\lambda(x)=
 $$
 $$\lim_{\epsilon \to  O+}\int_{\Delta_{\epsilon}}e^{\sum_{k=1}^{\infty}\frac{1}{2^k \epsilon}} f(x)
d\lambda(x)= \lim_{\epsilon \to  O+}\frac{1}{\lambda(\Delta_{\epsilon})}\int_{\Delta_{\epsilon}}
f(x) d\lambda(x).
$$
By Lemma 5.3 we know that
$$
\lim_{\epsilon \to  O+}\frac{1}{\lambda(\Delta_{\epsilon})}\int_{\Delta_{\epsilon}}
f(x) d\lambda(x)=f({\bf 0}).
$$
For $\alpha, \beta \in R$ and $f,g \in
~{C}(R^{\infty})$, we have
$$
\delta(\alpha f+\beta g)=(\delta)\int_{R^{\infty}}\delta(x)
(\alpha f(x)+\beta g(x)) d\lambda(x)= \lim_{\epsilon \to
 O+}\int_{R^{\infty}}\eta_{\epsilon}(x) (\alpha f(x)+\beta g(x))
d\lambda(x)=
$$
$$
\alpha \lim_{\epsilon \to
 O+}\int_{R^{\infty}}\eta_{\epsilon}(x)f(x) d\lambda(x)+\beta
\lim_{\epsilon \to 0+}\int_{R^{\infty}}\eta_{\epsilon}(x) g(x)
d\lambda(x)=\alpha f({\bf 0})+\beta g({\bf 0})=\alpha
\delta(f)+\beta \delta(g).
$$
This ends  the proof of the theorem.
~

Distributions are a class of linear functionals that map a set of
all test functions (conventional and well-behaved functions) onto
the set of real numbers. In the simplest case, the set of test
functions considered is $D(R^{\infty})$, which is the set of
smooth (infinitely differentiable) functions $\varphi : R^{\infty} \to
R$. Then, a distribution $d$ is a linear mapping $D(R^{\infty})\to
R$. Instead of writing $d(\varphi)$, where $\varphi$ is a test function in
$D(R^{\infty})$, it is conventional to write $\langle d,\varphi
\rangle$.

A simple example of a distribution is the Dirac delta functional
$\delta$, defined by $$\delta(\varphi) = \left\langle \delta,
\varphi \right\rangle = \varphi(0).$$

We have proved that Dirac delta functional $\delta$ is given by
the Dirac delta integral as follows
$$\delta(\varphi)=(\delta)\int_{R^{\infty}}\delta(x) \varphi(x)
d\lambda(x).$$

 There are straightforward mappings from both
locally integrable functions and probability distributions to
corresponding distributions, as discussed below. However, not all
distributions can be formed in this manner.

Suppose that $f : R^{\infty} \to R$ is a locally integrable
function, and let $\phi : R^{\infty} \to R$ be a test function in
$D(R^{\infty})$. We can then define a corresponding distribution
$T_f$ by
$$\left\langle T_{f}, \varphi \right\rangle =
\int_\mathbf{R^{\infty}} f(x) \varphi(x) \lambda (x).
$$
This integral is a real number which depends linearly and
continuously on $f$. This suggests the requirement that a
distribution should be a linear and continuous functional on the
space of test functions $D(R^{\infty})$, which completes the
definition. In a conventional abuse of notation, $f$ may be used
to represent both the original function $f$ and the distribution
$T_f$ derived from it. Similarly, if $\mu$ is a Radon measure on
$R^{\infty}$ and $f$ is a test function, then a corresponding
distribution $T_{\mu}$ may be defined by
$$\left\langle T_\mu, \varphi
\right\rangle = \int_{\mathbf{R^{\infty}}} \varphi\, d \mu.
$$

This integral depends continuously and linearly on $\varphi$, so
that $T_{\mu}$ is a distribution. If $\mu$ is an absolutely
continuous measure with respect to Baker measure $\lambda$ with
density $f$, then this definition is the same as the one for
$T_f$, but if $\mu$ is not absolutely continuous it gives a
distribution that is not associated with a function. For example,
if $P$ is the point-mass measure on $R^{\infty}$ that assigns $P$
measure one to the singleton set ${0}$ and measure zero to sets
that do not contain zero, then
$$
\int_{\mathbf{R^{\infty}}} \varphi \, dP =
(\delta)\int_{R^{\infty}}\delta(x) \varphi(x) d\lambda(x)=
\varphi(0),
$$ so
$T_P =\delta$ is the Dirac delta functional.

It is well known that the $n$-dimensional Dirac delta function
satisfies the following scaling property for a non-zero scalar
$\alpha$:
$$
(\delta) \int_{ R^n} \delta(\alpha x)\,dx =|\alpha|^{-n}
$$
 and so
$$
\delta(\alpha x) =|\alpha|^{-n}\delta(x).
$$

We have the direct generalization of that property in the case of
infinite dimension.

\begin{thm} The infinite dimensional Dirac delta function
satisfies the following scaling property for a non-zero scalar
$\alpha$:
$$(\delta)\int_{R^{\infty}}\delta(\alpha x)
d\lambda(x)=|\alpha|^{- \infty}.
$$
\end{thm}
{\bf Proof.}
We have
$$
(\delta)\int_{R^{\infty}} \delta(\alpha x) d\lambda(x)=
\lim_{\epsilon \to  O+}\int_{R^{\infty}}\eta_{\epsilon}(\alpha
x) d\lambda(x)=
$$
$$\lim_{\epsilon \to  O+}\int_{R^{\infty}}
\big[e^{\sum_{k=1}^{\infty}\frac{1}{2^k \epsilon}} \times \chi_{\Delta_{\epsilon}}(\alpha
x)+0 \times \chi_{R^{\infty} \setminus \Delta_{\epsilon}}(\alpha
x)\big]
 d\lambda(x)=
 $$
 $$\lim_{\epsilon \to  O+}e^{\sum_{k=1}^{\infty}\frac{1}{2^k \epsilon}} \Delta_{\epsilon}
e^{\frac{1}{\epsilon}} d\lambda(x)=
$$
$$
lim_{\epsilon \to
 O+}\frac{1}{\lambda(\Delta_{\epsilon})}\int_{\frac{1}{\alpha}\Delta_{\epsilon}}
d\lambda(x)= lim_{\epsilon \to
 O+}\frac{\lambda(\frac{1}{\alpha}\Delta_{\epsilon})}{\lambda(\Delta_{\epsilon})}.
$$
Notice that $\lambda(\frac{1}{\alpha}\Delta_{\epsilon})=0$ if
$|\alpha|> 1$, $=e^{-\sum_{k=1}^{\infty}\frac{1}{2^k \epsilon}}$ if $|\alpha|=1$ and
$=+\infty$ if $|\alpha|<1$.

Hence, the latter equality can be rewrited as follows
$$(\delta)\int_{R^{\infty}}\delta(\alpha x)
d\lambda(x)=|\alpha|^{- \infty}.
$$
This ends the proof of the theorem.

\begin{thm} The infinite dimensional Dirac delta function is an even distribution, in
the sense that
$$(\delta)\int_{R^{\infty}}\delta(-x) f(x)
d\lambda(x)=(\delta)\int_{R^{\infty}}\delta(x) f(x)
d\lambda(x)$$ for $f \in ~{C}(R^{\infty}),$ which
is homogeneous of degree $-1$.
\end{thm}

The validity of Theorem 5.6 follows from the fact asserted that
$-\Delta_{\epsilon}=\Delta_{\epsilon}$ for $\epsilon >0$ and the
invariance of $\lambda$ with respect to a transformation $T :
R^{\infty} \to R^{\infty}$ defined by $T(x)=-x$.

\begin{thm} ( sifting property ) The following equality
$$(\delta)\int_{R^{\infty}}\delta(x-T) f(x)
d\lambda(x)=f(T)$$holds for $f \in
~{C}(R^{\infty}).$
\end{thm}
{\bf Proof.}
By virtue of Lemma 5.3, we have
$$(\delta)\int_{R^{\infty}}\delta(x-T) f(x)
d\lambda(x)= \lim_{\epsilon \to
 O+}\int_{R^{\infty}}\eta_{\epsilon}(x-T) f(x) d\lambda(x)= $$
$$\lim_{\epsilon \to  O+}\int_{R^{\infty}}
\big[e^{\sum_{k=1}^{\infty}\frac{1}{2^k \epsilon}} \times \chi_{\Delta_{\epsilon}}(x-T)+0
\times \chi_{R^{\infty} \setminus \Delta_{\epsilon}}(x-T)\big]
 f(x) d\lambda(x)=
 $$
$$\lim_{\epsilon \to  O+}\int_{R^{\infty}}
\big[e^{\sum_{k=1}^{\infty}\frac{1}{2^k \epsilon}} \times \chi_{\Delta_{\epsilon}+T}(x)+0
\times \chi_{R^{\infty} \setminus (\Delta_{\epsilon}+T)}(x)\big]
 f(x) d\lambda(x)=
 $$
$$\lim_{\epsilon \to  O+}\int_{R^{\infty}}
e^{\sum_{k=1}^{\infty}\frac{1}{2^k \epsilon}} \times \chi_{\Delta_{\epsilon}+T}(x)
 f(x) d\lambda(x)=
  \lim_{\epsilon \to  O+}e^{\sum_{k=1}^{\infty}\frac{1}{2^k \epsilon}}\int_{\Delta_{\epsilon}+T}
  f(x) d\lambda(x)=
 $$
$$ \lim_{\epsilon \to  O+}\frac{\int_{\Delta_{\epsilon}+T}
  f(x) d\lambda(x)}{\lambda(\Delta_{\epsilon}+T)}=f(T).
 $$

\begin{thm} For $\epsilon >0$, let $(Y_n(\epsilon))_{n \in N}$ be an increasing family of finite subsets of $\Delta_{\epsilon}$ which is uniformly distributed in the $\Delta_{\epsilon}$. Let $f \in \mathcal{C}(R^{\infty})$. Then the following formula
$$
\lim_{\epsilon \to  0+} \lim_{n \to \infty} \frac{\sum_{y \in Y_n(\epsilon)}f(y)}{\#(Y_n(\epsilon))}=f({\bf 0})
$$
holds true.
\end{thm}

{\bf Proof.}  By Theorem 3.12 we have

$$
\lim_{n \to \infty} \frac{\sum_{y \in Y_n(\epsilon)}f(y)}{\#(Y_n(\epsilon))}=
\frac{\int_{\Delta_{\epsilon}}f(x)d\lambda(x)}{\lambda(\Delta_{\epsilon})}.
$$
By Lemma 5.3 we get
$$
\lim_{\epsilon \to  0+} \frac{\int_{\Delta_{\epsilon}}f(x)d\lambda(x)}{\lambda(\Delta_{\epsilon})}=f({\bf 0}),
$$
which implies that
$$
\lim_{\epsilon \to  0+} \lim_{n \to \infty} \frac{\sum_{y \in Y_n(\epsilon)}f(y)}{\#(Y_n(\epsilon))}=f({\bf 0}).
$$

This ends the proof of the theorem.

~

\begin{cor} For $\epsilon >0$, let $(Y_n(\epsilon))_{n \in N}$ be an increasing family of finite subsets of $\Delta_{\epsilon}$ which is uniformly distributed in the $\Delta_{\epsilon}$. Let $\delta$ be  Dirac delta functional defined in $\mathcal{C}(R^{\infty})$. Then the following equality
$$
\delta(f) = \lim_{\epsilon \to {\bf 0}} \lim_{n \to \infty} \frac{\sum_{y \in Y_n(\epsilon)}f(y)}{\#(Y_n(\epsilon))}
$$
holds true  for each $f \in \mathcal{C}(R^{\infty})$.
\end{cor}

\medskip
\noindent {\bf Acknowledgment}. The representation of the Dirac delta function  in terms of the Baker measure can be extended  also in terms of an arbitrary  ordinary or standard infinite-dimensional Lebesgue measure in $R^{\infty}$.

\end{document}